\documentclass{article}
\usepackage{a4wide}
\usepackage{lmodern}
\usepackage[T1]{fontenc}
\usepackage[latin9]{inputenc}
\usepackage{amsmath}
\usepackage{color}
\usepackage{graphicx}
\usepackage{amssymb}
\usepackage{xspace}
\usepackage{hyperref}
\usepackage{bbm}
\usepackage{oupbib}
\usepackage{todonotes}
\usepackage[draft]{listlbls}

\newcommand{\halm}{\hspace*{\fill} $\Box$\par}
\newtheorem{expl}{Example}[section]
\newenvironment{ex}{\begin{expl}\rm}{\halm\end{expl}}

\newcommand{\exref}[1]{\mbox{Example~\ref{ex:#1}}}
\renewcommand{\d}{\mbox{\rm d}}

 \newcommand{\cX}{{\cal
    X}} \newcommand{\cE}{{\cal E}} 
\newcommand{\cT}{{\cal T}} \newcommand{\cP}{{\cal P}}
 
 \newcommand{\R}{\mathbb{R}}

\newcommand{\bbl}{\mbox{$\mathbbm 1$}}

\newcommand{\cip}{\mbox{$\perp\!\!\!\perp$}}
\newcommand{\indo}[2]{\mbox{$#1 \,\cip\, #2$}}

\newcommand{\reals}{\R} 
\newcommand{\transp}{^{\rm T}} 
\newcommand{\itref}[1]{\mbox{\ref{it:#1}}}
 
 \newcommand{\xbar}{\mbox{$\overline x$}}

\newcommand{\wbar}{\mbox{$\overline w$}}
\newcommand{\ebar}{\mbox{$\overline e$}}
\newcommand{\Ebar}{\mbox{$\overline E$}}

\newcommand{\ie}{{\em i.e.\/}\xspace}
\newcommand{\eg}{{\em e.g.\/}\xspace}

\renewcommand{\eqref}[1]{\mbox{(\ref{eq:#1})}}
\newcommand{\secref}[1]{\mbox{\S$\,$\ref{sec:#1}}\xspace}
\newcommand{\Secref}[1]{\mbox{Section~\ref{sec:#1}}}

\title{Fiducial inference then and now}

\author{Philip Dawid\thanks{University of Cambridge}}

\date{\today}

\begin{document}

\maketitle
\begin{center}
  Dedicated to the memory of Mervyn Stone, 1932--2020
\end{center}

\begin{abstract}
  \noindent We conduct a review of the fiducial approach to
  statistical inference, following its journey from its initiation by
  R. A. Fisher, through various problems and criticisms, on to its
  general neglect, and then to its more recent resurgence.  Emphasis
  is laid on the functional model formulation, which helps clarify the
  very limited conditions under which fiducial inference can be
  conducted in an unambiguous and self-consistent way.\\

\end{abstract}

\section{Introduction}
\label{sec:intro}
According to \textcite{zabell}, ``the fiducial argument stands as
Fisher's one great failure'', a sentiment that has been echoed by
others.  Fisher never constructed a fully-fledged theory of fiducial
inference, but developed his ideas by means of examples and {\em ad
  hoc\/} responses to increasingly complex problems or challenges
raised by others.  Few other statisticians have taken the fiducial
argument seriously, and after some sporadic activity (mostly critical)
in the decades following Fishere's introduction of the idea in 1930,
it almost completely disappeared from the scene.  The trio of
Encyclopedia articles \textcite{awfe:ess,rjb:ess,ms:ess} is a useful
resource for the state of the enterprise up to 1982.
% ; another interesting read is \textcite{zabell}.
In recent times, however, there has been a resurgence of interest in
the fiducial programme, as evidence by works such as \textcite{hannig}
and \textcite{martin/liu:book}, and the success of the series of
annual {\em Bayesian, Fiducial \& Frequentist\/} (BFF) conferences,
since 2014.

In this article I give a personal review of the main contributions,
positive and negative, to fiducial inference, in both earlier and
later periods.  In \secref{fisher}I describe the original argument,
centered on inference for a correlation coefficient.  \Secref{pivot}
introduces an extension to more complex problems, based on the idea of
a pivotal function.  Such a function arises naturally when the problem
possesses properties of invariance under a group of transformation, as
described in \secref{group}.  A variant of this is Fraser's structural
model described in \secref{struct}, while a further extension is the
functional model of \secref{func}, which forms a basis for the rest of
the article.  In \secref{margin} we show that, under certain
conditions, two different routes to marginalizing a fiducial
distribution give the same answer.

In \secref{diff} we start to see some problems with the fiducial
argument.  In particular, when attempting to condition in a fiducial
distribution, we again have two possible routes, but they generally
yield different answers.

To this point we have only considered simple models, essentially those
where the dimensions of the parameter and the data are the same.
Non-simple models require additional conditioning, as described in
\secref{nonsimp}.  However this requires an additional property,
partitionability, in the absence of which there is no well-defined
fiducial distribution.

\Secref{nifm} considers cases in which the fiducial argument fails to
yield a distribution for the parameter, but only a distribution for a
set containing the parameter.  This is linked to the Dempster-Shafer
theory of belief functions, and has been a focus of recent work.  In
this case too a partitionability property is required for well-defined
inference.

Some concluding thoughts are gathered in \secref{conc}.

\section{Fisher's original fiducial argument}
\label{sec:fisher}

The fiducial argument was introduced by \textcite{fisher:invprob} by
means of the following example.

We have $n$ observations from a bivariate normal distribution.  Let
random variable $R$ be the sample correlation, and let parameter
variable $\Phi$ be the population correlation.  Then the sampling
distribution of $R$ depends only on ($n$ and) the value $\phi$ of
$\Phi$.  The form of this distribution \cite{fisher:cor} is not
expressible by means of simple functions; however the rest of Fisher's
argument does not involve the specific form of this distribution.
Indeed, letting $F(r;\phi) = \Pr(R \leq r \mid \Phi = \phi)$ be the
cumulative distribution function of $R$ when $\Phi=\phi$, the general
argument applies to any problem satisfying the following sufficient
(but not entirely necessary) regularity conditions:
\renewcommand{\theenumi}{(\roman{enumi})}
\begin{enumerate}
\item \label{it:domain} The range of each of $R$ and $\Phi$ is an
  open interval in the real line
\item \label{it:conts}
  $F(r;\phi)$ is a continuous function of each of its arguments
\item \label{it:rup} For fixed $\phi$, $F(r;\phi)$ strictly increases
  as $r$ increases, taking all values in $(0,1)$
\item \label{it:phidown} For fixed $r$, $F(r;\phi)$ strictly decreases
  as $\phi$ increases, taking all values in $(0,1)$
\end{enumerate}

It follows from the probability integral transformation \cite{angus}
that, for all $\phi$, the distribution of $F(R;\phi)$, given
$\Phi = \phi$, is uniform on $[0,1]$.  That is,
$E = F(R;\Phi) \sim U[0,1]$, independently of $\Phi$: $\indo E \Phi$.

Since, for any $\gamma\in[0,1]$, when $\Phi=\phi$
\begin{equation}
  \label{eq:cdf}
 \Pr\{F(R;\phi) \leq \gamma\} = \gamma,
\end{equation}
a level-$\gamma$ confidence set for $\Phi$ is, for observed $R=r$,
$I(r;\gamma) := \{\phi: F(r;\phi) \leq \gamma\}$.  In fact, because
$F(r;\phi)$ is a decreasing function of $\phi$, this is an upper
confidence interval: $I(r;\gamma) = [\phi(\gamma),\infty)$, where
$F(r;\phi(\gamma)) = \gamma$.

Fisher now takes this argument further.  He regards the uniform
distribution for $E = F(R;\Phi)$, and so \eqref{cdf}, as remaining
valid, even after observing data wih $R=r$.  Equivalently, he takes
$\indo E R$ (compare the sampling property $\indo E \Phi$).
Thus he assumes, for any $r$, 
\begin{equation}
\label{eq:uniffid}
F(r,\Phi) \sim U[0,1].
\end{equation}
This argument has now assigned to the parameter $\Phi$ the status of a
random variable.  Indeed, after observing $R=r$, \eqref{uniffid}, in
conjunction with \itref{conts} and \itref{phidown}, implies
\begin{equation}
  \label{eq:phis}
  \Pr(\Phi\leq \phi) = \Pr\{F(r;\Phi)\geq F(r;\phi)\} = 1-F(r;\phi).
\end{equation}
On account of \itref{phidown} this yields a full ``fiducial
distribution function'' for $\Phi$; under differentiability, the
associated ``fiducial density'' of $\Phi$ is
$ -\frac{\partial}{\partial \phi} F(r;\phi)$.  In particular,
$\Pr\{\Phi \in I(r;\gamma)\} = \Pr\{\Phi \geq
\phi(\gamma)\}=F\{r,\phi(\gamma)\} = \gamma$, so transforming a
confidence statement for $\Phi$ into a probability statement for
$\Phi$---an interpretation of a confidence interval that is typically
castigated as a showing a gross misunderstanding of its nature.

Fisher, noting that \eqref{phis} does not follow from standard
probability arguments, termed it a ``fiducial probability''.  However,
he appeared to believe that (subject to some caveats---see
\secref{val} below) it can still be interpreted as a regular
probability.  \textcite{savage:eggs} memorably described the fiducial
argument as ``a bold attempt to make the Bayesian omelet without
breaking the Bayesian eggs''.  \textcite{dvl:fidbayes} considered the
general case \itref{domain}--\itref{phidown} of Fisher's construction,
and showed that the fiducial distribution does not arise as a Bayesian
posterior, for any prior, except in the special case of a location
model with uniform prior on the location parameter.

\section{Pivotal inference}
\label{sec:pivot}
A more general fiducial construction relies on the existence of a
pivot \cite{barnard}.

Let $X$ have distribution governed by parameter $\Theta$.  A {\em
  pivot\/} $E$ is a function of $X$ and $\Theta$ with known
distribution $P_0$, not depending on the value $\theta$ of $\Theta$:
$\indo E \Theta$, with $E \sim P_0$.

\begin{ex}
  \label{ex:cdfpiv}
  In \secref{fisher}, $E = F(R;\Phi)$ is a pivot, with distribution
  $U[0,1]$.  We have seen how this can be used to supply a fiducial
  distribution.
\end{ex}

\begin{ex}
  \label{ex:normpivs}
  For a sample of size $n$ from the normal distribution with mean $M$
  and variance $\Sigma^2$, having sample mean $\overline X$ and sample
  variance $S^2$, $E_1= (\overline X - M)/\Sigma$ is a pivot, wih the
  normal distribution $N(0,1/n)$; $E_2 = S/\Sigma$ is a pivot, with
  distribution $\sqrt{\chi^2_{n-1}/(n-1)}$; and
  $E_3 =\sqrt{n} E_1/E_2 = (\overline X - M)/(S/\sqrt n)$ is a pivot,
  with the Student distribution $t_{n-1}$.
\end{ex}
Given a suitable pivot $E = f(X,\Theta)$, a fiducial distribution is
obtained by regarding the distribution $P_0$ of $E$ as still relevant,
even after observing the data.  Equivalently, instead of
$\indo E \Theta$, we regard $\indo E X$.  Thus after observing
$X = x$, we suppose $f(x,\Theta) \sim P_0$.  When $\Theta$ (or a desired
function $\Psi$ of $\Theta$) can be expressed as a function of
$f(x,\Theta)$---the case of {\em invertibility\/}---this delivers a
fiducial distribution for $\Theta$ (or $\Psi$.)

In \exref{normpivs}, using $E_3$ we have
$(\overline x - M)/(s/\sqrt n)\sim t_{n-1}$, which can be solved as
$M = \overline x - (s/\sqrt n)E_3$, yielding fiducial distribution
$M \sim \overline x + (s/\sqrt n) t_{n-1}$.  Similarly using $E_2$ we
obtain a fiducial distribution for $\Sigma$:
$\Sigma = s/E_2\sim s/\sqrt{\chi^2_{n-1}/(n-1)}$.  The pivot $E_1$, by
itself, is not invertible, and does not yield a fiducial distribution.
However, the bivariate pivot $(E_1, E_2)$ is invertible, and yields a
joint fiducial distribution for $(M,\Sigma)$, represented by
\begin{eqnarray}
  \label{eq:mrep}
  M &=& \overline x -  s{E_1}/{E_2}\\
  \label{eq:srep}
  \Sigma &=& s/{E_2}
\end{eqnarray}
where still $E_1\sim N(0,1/n)$, $E_2 \sim \sqrt{\chi^2_{n-1}/(n-1)}$
(these moreover retaining their sampling distribution independence).
In particular, the induced marginals for $M$ and $\Sigma$ agree with
those above based directly on $E_3$ and $E_2$, as above.  However, $M$
and $\Sigma$ are not independent.  Learning $\Sigma = \sigma$ is
equivalent to learning $E_2 = s/\sigma$.  This does not change the
$N(0,1/n)$ distribution of $E_1$, and we now have
$M = \overline x - \sigma E_1$.  So we have conditional fiducial
distribution
\begin{equation}
  \label{eq:mmids}
  M \mid (\Sigma=\sigma) \sim N(\overline x, \sigma^2/n).  
\end{equation}

\subsection{Validity}
\label{sec:val}
There will typically be many available pivotal functions.  For
instance, in \exref{normpivs} we could retain just the first $n/2$
(say) observations, and use the sample mean and variance computed from
these.  In his early writings, Fisher insisted that, to make use of
all the available information, a fiducial distribution should be based
on the minimal sufficient statistic.

In addition, Fisher indicated that a fiducial distribution should be
regarded as yielding an appropriate inference only if the following
vaguely stated conditions are satisfied:
\begin{enumerate}
\item \label{it:noprior} there is no available prior information about
  the unknown parameter;
\item \label{it:uninf} ({\em ``principle of irrelevance''\/}): the
  data are uninformative about the pivot.
\end{enumerate}
While \itref{noprior} obviously precludes having a Bayesian prior
distribution, its intended scope is much wider (and much vaguer).  As
for \itref{uninf}, \textcite{hacking:book} (see also
\textcite{harris}) attempted to rigorize it as requiring that the
likelihood function based on any data, re-expressed as function of the
data and the pivot, be the same (up to proportionality) for any data.
He then showed that, in the univariate case, this holds if and only if
(possibly after transformation) we have a location model, in which
case, as observed by \textcite{dvl:fidbayes}, the fiducial
distribution agrees with the formal posterior based on an improper
uniform prior.

\section{Group-structured models}
\label{sec:group}
Let observable $X$ take values in a space $\cX$, and identifiable
parameter $\Theta$ take values in $\cT$.  We suppose $\cT$ can be
identified with a group $G$ of transformations acting on $\cX$.  We
denote the image of $x\in\cX$ under $g\in G$ by $g\circ x$, and
further suppose the group action is {\em exact\/}, so that, given
$x_0, x_1\in \cX$ there is at most one $g\in G$ such that
$g\circ x_1 = x_2$ (this condition can be relaxed: see
\textcite{bondar}).

We shall investigate cases in which the family
$\cP = \{P_\theta:\theta\in\cT\}$ of distributions over $\cX$ is {\em
  equivariant\/} under the action of $G$: that is, if
$X\sim P_\theta$, then $g\circ X\sim P_{g\theta}$ (where $g\theta$ is
the group product).

\subsection{Simple group model}
\label{sec:simplest}
In the simplest case, $G$ acts transitively (as well as exactly) on
$\cX$, so that, for any $x_1,x_2 \in \cX$ there exists exactly one
$g\in G$ such that $x_2 = g\circ x_1$.  Fix some $x_0\in \cX$, and
henceforth identify any $x\in \cX$ with the unique $g\in G$ such that
$x=g\circ x_0$.  We can thus take $\cX = G$, with $g\circ x$ becoming
the group product $gx$.  Let $P_0 = P_\iota$, with $\iota$ the
identity element of $G$.  Define $E := \Theta^{-1} X$.  Then, by
equivariance, conditional on $\Theta=\theta$, $E \sim P_0$.  Hence $E$
is a pivot, and could be used to construct a fiducial distribution:
after observing $X=x$, take $\Theta^{-1}x \sim P_0$.  It can be shown
that this construction satisfies Hacking's version of the principle of
irrelevance.  Moreover \cite{fraser:1961}, assuming $G$ is locally
compact, the resulting fiducial distribution is identical to a
Bayesian posterior distribution, based on the (typically improper)
right-invariant distribution (right Haar measure) for $\Theta$ over
$G$.

\begin{ex}
  \label{ex:normpivs2}
In \exref{normpivs}, we can consider both $(\overline X, S)$ and
$(M,\Sigma)$ as elements of the location-scale group, with
multiplication $(a,b)(A,B) = (a+bA,bB)$, and identity $\iota = (0,1)$.
% Then
% $(\mu,\sigma)(\overline x, s) = (\mu + \sigma \overline x, \sigma
% s)$.
Under $P_0$, $\overline X\sim N(0,n^{-1})$ and
$S\sim \sqrt{\chi^2_{n-1}/(n-1)}$, independently; then
$(\mu,\sigma)(\overline X, S) = (\mu+\sigma\overline X, \sigma S) \sim
(N(\mu, \sigma^2/n), \sigma\sqrt{\chi^2_{n-1}/(n-1)})$
(independently), which is $P_{(\mu,\sigma)}$.  Thus we have
equivariance.

We have pivot
\begin{eqnarray*}
  E = (E_1,E_2) &=& (M,\Sigma)^{-1}(\overline X, S)\\
                &=& \left(\frac{\overline X-M}{\Sigma}, \frac S {\Sigma}\right).
\end{eqnarray*}
We recover the same joint fiducial distribution represented by
\eqref{mrep} and \eqref{srep}.  Moreover, this is the same as the
posterior distribution based on the right-invariant prior, having
density element $\d\mu\,\d\sigma/\sigma$.
\end{ex}

\begin{ex}
  \label{ex:prog} Let $G$ be the group of lower triangular matrices
  with positive diagonal.  An observable random $2\times 2$ matrix $S$
  has the Wishart distribution $W(\nu;\Sigma)$, where $\nu\geq 2$ and
  $\Sigma$ is positive definite.  Then $S$ is almost surely
  non-singular.  We can alternatively represent $S$ by the unique
  $L\in G$ such that $S = LL\transp$, and similarly $\Sigma$ by
  $\Lambda\in G$ with $\Sigma = \Lambda\Lambda\transp$.  We write the
  implied distribution of $L$, depending on $\Lambda$, as
  $L \sim {\cal L}(\nu;\Lambda)$.  It is then easy to see that, for
  fixed $A\in G$, $AL \sim {\cal L}(\nu;A\Lambda)$, so that the
  problem is equivariant under $G$.  It follows that a pivot is
  $E = \Lambda^{-1}L$, with distribution ${\cal L}(\nu;I)$---under
  which \cite{mauldon} the non-zero entries of $E$ are independent,
  with $E_{11}\sim\sqrt{\chi^2_\nu}$,
  $E_{22} \sim\sqrt{\chi^2_{\nu-1}}$, and $E_{21}\sim N(0,1)$.  So the
  fiducial distribution of $\Sigma = \Lambda\Lambda\transp$, given
  data $S =s = ll\transp$, is that of $l(E\transp E)^{-1}l\transp$
  when $E\sim{\cal L}(\nu;I)$.  In this case the right-invariant prior
  density element, under the action of $G$, can be expressed in terms
  of the entries of $\Sigma$ as
  \begin{equation}
    \label{eq:haar}
  (\d\sigma_{11}/\sigma_{11})\,\d\sigma_{12}\,\d\sigma_{22},  
\end{equation}
and the fiducial distribution of $\Sigma$ agrees with its posterior,
based on \eqref{haar} as prior.
\end{ex}

\subsection{Structural models}
\label{sec:struct}
\textcite{fraser:1961,fraser:book} has a somewhat different take on
group-structured models, which
% again avoids the non-uniqueness problem of \secref{mauldon} by
% taking
takes the group structure as part of the specification of the problem.
He posits, as part of the very set-up, a nominated group $G$ acting on
$\cX$, and an ``error variable'' $E$, with known distribution $P_0$
over $\cX$.  Both $E$ and the parameter $\Theta$, which takes values
in $G$, are regarded as having independent existence.  The observable
$X$ is then {\em defined\/} by $X=\Theta \circ E$.  There is thus
additional algebraic structure, over and above the implied parametric
family of distributions for $X$ given $\Theta=\theta$.  This extended
structure is termed a ``structural model''.  

Since the implied distributional model is equivariant under $G$, we
can now construct a fiducial distribution as in \secref{simplest} (see
also \secref{nontrans} below for the non-transitive case).  Fraser
terms this a ``structural distribution''.

Note that, as demonstrated in \secref{mauldon} below, distinct
structural models can correspond to the same distributional model, and
yield different structural distributions for its parameter.  Since a
structural model is considered to comprise more than just its induced
distributional model---including, in particular, specification of the
group $G$ as a key ingredient---this is not regarded as an
inconsistency.

\section{Functional models}% (Fraser; Dawid \& Stone)}
\label{sec:func}

\textcite{apd/ms:funcmod} propose the {\em functional model\/}, a
generalization of the structural model.  We have arbitrary sample
space $\cX$ and parameter space $\Theta$.  We again consider an
``error variable'' $E$, taking values in a space $\cE$ that now may be
different from $\cX$.  The observable $X$ is defined, algebraically,
as $X = f(\Theta,E)$, where $f:\cT\times\cE\rightarrow\cX$ is a
specified function, and $E$ has a known distribution $P_0$ over $\cE$,
independently of the value of $\Theta$.  For simplicity we denote the
function simply by $X = \Theta\circ E$.  Then the distribution
$P_\theta$ of $X$ given $\Theta=\theta$ is that of $\theta\circ E$
where $E\sim P_0$.

The structural approach, which is a special case of the functional
approach, identifies $\theta\in\cT = G$ with the function (an element
of $G$) $e \mapsto \theta\circ e$ on $\cX$.  In the functional
approach, by contrast, it is more helpful to consider $e\in \cE$ as
the function $\theta\mapsto \theta\circ e$, mapping $\cT$ into $\cX$.
(Note that, as a function, $e$ is written to the right of its argument
$\theta$.)

\subsection{Simple functional model (SFM)}
\label{sec:sfm}
In the simplest case, for any $x\in\cX$, $e\in\cE$, there exists
exactly one $\theta$ such that $x = \theta\circ e$: we write
$\theta = x \circ e^{-1}$, since this determines the inverse function
$e^{-1}:\cX\rightarrow \cT$.  In this case the fiducial distribution,
for data $x$, is obtained from $\Theta = x\circ E^{-1}$, with
$E \sim P_0$.

In the special case that $E$ can be expressed as a function of
$(X,\Theta)$, it serves as a pivot.  The model is then termed {\em
  pivotal\/}, and the fiducial distribution agrees with that
constructed as in \secref{pivot}.

\subsubsection{Monotonic functional model}
\label{sec:decreasing}
When $\cX = \cT = \reals$ and each $e$ acts as a strictly monotonic
function of $\theta$, the fiducial distribution is fully determined by
the distributional model for $X$ given $\Theta$: in particular, the
finer details of the functional model do not enter.  Thus when $e$ is
a decreasing function, the fiducial probability
$\Pr(\Theta \leq \theta) = P_0(x\circ E^{-1 }\leq\theta) = P_0(x \geq
\theta\circ E) = \Pr_\theta(X \leq x)$.  One can show that, under
regularity conditions parallel to those in \secref{fisher}, and by a
similar argument, a 1-sided fiducial interval is also a confidence
interval.

\begin{ex}
  \label{ex:nonpivot}
  Let $\cX = \cT = \reals$, $\cE = (\reals^+)^3$.  The function
  $x = \theta\circ e$ is given by $x = (\theta e_1 + e_3)/e_2$, which
  is strictly increasing in $\theta$.  This model is not pivotal, but
  we can solve for $\theta$:
  $\theta = x\circ e^{-1} = (xe_2 - e_3)/e_1$.  So the fiducial
  distribution of $\Theta$, for data $X=x$, is that of
  $(xE_2 - E_3)/E_1$, with $E$ having its initially assigned
  distribution $P_0$.

  As a special case, suppose that, under $P_0$,
  $E_1\sim \sqrt{\chi^2_{n-1}}$, $E_2\sim \sqrt{\chi^2_{n-2}}$, and
  $E_3\sim N(0,1)$, all independently.  Define $R= X/\sqrt{1+X^2}$,
  $\Phi:= \Theta/\sqrt{1+\Theta^2}$.  It then turns out (see
  \exref{prog2} below) that (compare \secref{fisher}):
  \begin{enumerate}
  \item \label{it:corsamp} the sampling distribution of $R$ is that of
    a sample correlation coefficient, based on $n$ independent
    observations from a bivariate normal distribution with population
    correlation coefficient $\Phi$
  \item \label{it:corfid} the fiducial distribution of $\Phi$ agrees
    with Fisher's fiducial distribution, based on $R$.  This follows
    from \itref{corsamp} as a consequence of the monotonic structure
    of this model.
  \end{enumerate}
\end{ex}

\section{Marginalization consistency}
\label{sec:margin}
In a SFM $X = \Theta\circ E$, with $E\sim P_0$, let $W = w(X)$ be a
function of $X$.  Suppose that $w=w(\theta\circ e)$ can be expressed
as a function of $\omega$ and $e$, where $\omega = \omega(\theta)$ is
some function of $\theta$: we write this function as $w = \omega*e$.
With $\Omega := \omega(\Theta)$ we thus have a new model
$W = \Omega*E$: in particular, the sampling distribution of $W$
depends only on the value of $\Omega$.  We require that this model
itself be a SFM, so that, given $(w,e)$, we can solve $w=\omega*e$ for
$\omega$, which solution we write as $\omega = w*e^{-1}$ (though the
function $e^{-1}$, now acting on $\omega$, has a different meaning
here than in $\theta = x\circ e^{-1}$.)  We term the model
$W = \Omega*E$ a {\em reduction\/} of
$X = \Theta\circ E$.\footnote{Appendix~A2 of \textcite{apd/ms:funcmod}
  characterises such a reduction in terms of group actions.  In
  particular, if the initial SFM is structural, with $\cT$ a group $G$
  of transformations of $\cX=\cE$, a reduction is obtained by taking
  $W$ and $\Omega$ as maximal invariants under a subgroup $K$ of $G$,
  acting on $\cX$ and $\cT$ respectively.}

Given data $X=x$, we have two different routes to computing the
fiducial distribution of $\Omega = \omega(\Theta)$:
\begin{enumerate}
\item \label{it:top} Obtain the fiducial distribution of $\Theta$
  based on data $X=x$, using the full SFM $X = \Theta\circ E$; then
  marginalize this to get the implied distribution of
  $\Omega = \omega(\Theta)$.
\item \label{it:bottom} Start from the reduced SFM $W = \Omega*E$, and
  obtain the associated fiducial distribution of $\Omega$, based on
  the reduced data $W=w(x)$.
\end{enumerate}
To see that these give the same result we argue as follows.  Route
\itref{top} represents $\Theta = x \circ E^{-1}$, and so produces the
distribution of $\omega(x\circ E^{-1})$.  Route \itref{bottom}
represents $\Omega = w(x) *E^{-1}$.  In both cases $E\sim P_0$.  Now
if $x = \theta\circ e$ then $\theta = x \circ e^{-1}$.  Also
$w(x) = \omega(\theta)*e$, so $\omega(\theta) = w(x)*e^{-1}$.  Hence
$\omega(x\circ e^{-1})= w(x) * e^{-1}$ whence
$\omega(x\circ E^{-1})= w(x) * E^{-1}$, showing that both routes yield
the same representation, and hence the same fiducial distribution, for
$\Omega$ (in particular, the marginal fiducial distribution of
$\Omega$ in the route~\itref{top} analysis must depend on the data $x$
only through $w = w(x)$).

\begin{ex}
  \label{ex:normpivs3}
  \exref{normpivs2} can be regarded as a SFM (in fact a structural
  model): $(\overline X, S) = (M,\Sigma)\circ(E_1,E_2)$, where $\circ$
  is group product in the location-scale group.  That is,
  \begin{eqnarray}
      \label{eq:xx}
    \overline X &=& M + \Sigma E_1\\
    \label{eq:ss}
    S &=& \Sigma\, E_2
  \end{eqnarray}
  with $E_1\sim N(0,1/n)$, $E_2\sim\sqrt{\chi^2_{n-1}/(n-1)}$,
  independently.
    
  Define $W = \overline X/S$, $\Omega = M/\Sigma$.  Then
  \begin{equation}
    \label{eq:wfm}
    W = \frac{\Omega + E_1}{E_2}.
  \end{equation}
  This is itself a SFM $W = \Omega*F$ (though not structural), so is
  a reduction of the original SFM.

  \textcite{demp:incons} noted that (as indeed follows from the above)
  the sampling distribution of $W$ depends only on $\Omega$, so that a
  fiducial distribution for $\Omega$ can be constructed from these
  univariate sampling distributions using Fisher's approach of
  inverting the distribution function.  By monotonicity, the
  route~\itref{bottom} analysis of the reduced model $W = \Omega*F$
  will also deliver this fiducial distribution.
  \textcite{demp:incons} further showed that this agrees with the
  distribution of $\Omega = M/\Sigma$ obtained by marginalising the
  joint fiducial distribution of $(M,\Sigma)$ represented by
  \eqref{mrep} and \eqref{srep}---as also arises from the
  route~\itref{top} analysis of the initial SFM
  \eqref{xx}--\eqref{ss}.  Here we see marginal consistency in action.
  
\end{ex}

\begin{ex}
  \label{ex:prog2}
  \exref{prog} can be regarded as the (structural) SFM $L=\Lambda E$
  (all lower triangular matrices), with $E_{11}\sim\sqrt{\chi^2_\nu}$,
  $E_{22} \sim\sqrt{\chi^2_{\nu-1}}$, and $E_{21}\sim N(0,1)$,
  independently.  Thus
  \begin{eqnarray*}
    L_{11} &=& \Lambda_{11}E_{11}\\
    L_{12} &=& \Lambda_{12}E_{11}+\Lambda_{22}E_{22}\\
    L_{22} &=& \Lambda_{22}E_{22}.
  \end{eqnarray*}
  Defining $X=L_{12}/L_{11}$, $\Theta = \Lambda_{12}/\Lambda_{11}$, we
  obtain a reduction $X = \Theta*E$, given by the SFM
  $X = (\Theta E_{11} + E_{12})/E_{22}$.  Note that (with minor
  notational changes) this is identical with the special case
  considered in \exref{nonpivot}.

  From \exref{prog} we have $S = LL\transp$, \ie
  $$\left(
    \begin{array}[c]{cc}
      S_{11} & S_{12}\\
      S_{12} & S_{22}
    \end{array}
  \right)
  = \left(
    \begin{array}[c]{cc}
      L_{11}^2 & L_{11}L_{12}\\
      L_{11}L_{12} & L_{12}^2 + L_{22}^2
    \end{array}
  \right),$$ with a similar expression for $\Sigma$ in terms of
  $\Lambda$.  The sample correlation based on $S$ is
  $$R := \frac{S_{12}}{\sqrt{S_{11}S_{22}}} = \frac X {\sqrt{1+X^2}},$$
  and similarly the population correlation is
  $\Phi := \Theta /\sqrt{1+\Theta^2}$.  The former identity explains
  the distribution of $R$ asserted in \itref{corsamp} of \exref{prog}.
  There we deduced that the fiducial distribution of $\Phi$, based on
  $R$, agrees with that derived by Fisher.  By marginalization
  consistency, this must also be true for the marginal distribution of
  $\Phi =\Sigma_{12}/\sqrt{\Sigma_{11}\Sigma_{22}}$ formed from the
  full fiducial distribution of $\Sigma$ given $S$, based on the
  lower-triangular structural model.

  \subsection{Marginalization paradox}
\label{sec:mp}
Although the marginalization consistency property seems to speak in
favour of fiducial inference, at least in some problems, it becomes a
problem for Bayesian inference with improper priors.  We know that the
full fiducial distribution, used in the route~\itref{top} analysis, is
also the Bayesian posterior, based on the right-invariant prior
distribution.  Thus the output of the route~\itref{top} analysis is
the marginal distribution of $\Omega$ in this Bayesian posterior.  By
marginalization consistency, the output of route~\itref{bottom}
analysis---which depends on the data $X$ only through $W$---must then
likewise agree with this marginal distribution.  It therefore seems
reasonable to believe that the marginal Bayesian distribution of
$\Omega$, depending as it does only on $W$, could arise as a Bayesian
posterior based on the likelihood from the reduced model for $W$
(depending only on $\Omega$).  But by the result of
\textcite{dvl:fidbayes} (see \secref{fisher}), if---as in both the
above examples---the reduced model has univariate $W$ and $\Omega$ but
is not equivalent to a location model, this can not be the case.  We
then have an example of a {\em marginalization paradox\/} \cite{dsz}
in improper Bayesian inference.

\section{Some difficulties}
\label{sec:diff}

\subsection{Choice of group}
\label{sec:mauldon}
\textcite{mauldon} pointed out a problem with \exref{prog}: if we
simply interchange the order in which we consider the variables
(equivalent to now using equivariance under the upper triangular,
rather than lower triangular, group), the analysis proceeds
essentially as before, but we obtain a different fiducial
distribution.  This can most easily be seen by noting that the
right-invariant prior \eqref{haar} is altered on interchanging the
suffices 1 and 2, leading to a different posterior, hence fiducial,
distribution for $\Sigma$.

One possible escape from this bind is not to allow the use of just any
group $G$ under which the statistical model happens to be equivariant,
but to specify an appropriate group as part of the very structure of
the problem---this thus requiring an additional ingredient in the
model, over and above its purely distributional properties.  This
tallies with the position adopted in Fraser's structural
modelling---see \secref{struct}.

Although use of the upper triangular group produces a different
fiducial distribution for $\Sigma$ than that based on the lower
triangular group, nevertheless, by an argument parallel to that of
\exref{prog2}, the implied distribution for $\Phi$ again agrees with
Fisher's, and thus is the same in both cases.
  \end{ex}

  \subsection{Marginalization inconsistency}
\label{sec:margincons}

\begin{ex}
    \label{ex:wilk}
    Consider the $n$-variate SFM $X = \Theta\circ E$ given by
  $X_i = \Theta_i + E_i$ ($i=1,\ldots,n$), with $E_i \sim N(0,1)$, all
  independently.  On observing $X=x$, the fiducial distribution has
  $\Theta_i \sim N(x_i,1)$, independently.

  Let $W = \sum_{i=1}^n X_i^2$, $\Omega = \sum_{i=1}^n \Theta_i^2$.
  The marginal fiducial distribution of $\Omega$, given data $x$,
  depends only on $w = \sum_{i=1}^n x_i^2 $; it is non-central
  $\chi^2$ with non-centrality parameter $w$: $\Omega\sim\chi^2_n(w)$.
  Also, the sampling distribution of $W$, when $\Theta = \theta$,
  depends only on $\omega = \sum_{i=1}^n \theta_i^2$; it is
  $\chi^2_n(\omega)$.  Nevertheless,
  $W = \sum_{i=1}^n (\Theta_i + E_i)^2$ can not be expressed as a
  function of $\Omega$ and $E$, so we do not have a reduction of the
  initial model.  We note that $W-\Theta$ has sampling expectation
  $n$, but fiducial expectation $-n$, which suggests a serious
  inadequacy in the marginalized fiducial distribution.

  We can attempt to derive a ``route~\itref{bottom}''-type fiducial
  distribution of $\Omega$, by Fisherian inversion of the distribution
  function of $W$ given $\Omega$.  In this case condition \itref{rup}
  of \secref{fisher} does not hold, and we obtain only an incomplete
  distribution---which obviously can not agree with the complete
  marginal fiducial distribution obtained from route~\itref{top}
  analysis, so we do not have marginalization consistency.  In fact,
  with $n=50$, we get 95\% central fiducial interval $(109,196)$ by
  marginalizing the full fiducial distribution to $\Omega$, compared
  to $(21,89)$ based on the distribution of $W$ given $\Omega$.
\end{ex}

This example indicates that sensible marginalization of a joint
fiducial distribution may not be possible when not based on a
reduction of a functional model.  \textcite{wilkinson} embraces
inconsistencies such as in this example by his {\em noncoherence
  principle\/}, which allows the overall joint fiducial to coexist
with the ``marginal'' based on the reduced data---which is not the
actual marginal.  But then fiducial distributions do not satisfy the
axioms of probabiity theory.

Other examples of marginalization inconsistency, evidenced by
incompatibilities between fiducial and confidence statements, are the
Behren-Fisher problem, looking at the difference between the means of
two normal distributions with different, unknown, variances, and the
Fieller-Creasy problem, looking at the ratio of two normal means, with
known variances \cite{wallace:bffc}.

\subsection{Conditional consistency?}
\label{sec:condcon}

\begin{ex}
  \label{ex:normpivs4}
  Consider again \exref{normpivs}, and suppose we want to construct a
  conditional fiducial distribution for $M$, given $\Sigma=\sigma$.
  Again we have two possible routes to do this:
  \begin{enumerate}    
  \item \label{it:cond1} Condition the joint fiducial distribution on
    $\Sigma = \sigma$, leading to \eqref{mmids}.
  \item Note than, when $\Sigma=\sigma$ is fixed, the sampling model
    now has
    \begin{eqnarray}
      \label{eq:aa1}
      \overline X &=& M +\sigma E_1\\
      \label{eq:aa2}
      S &=& \sigma E_2
    \end{eqnarray}
    On observing $(\overline x,s)$ we learn $E_2 = s/\sigma$, so
    should condition on this.  The same reasoning that led to
    \eqref{mmids} again applies, so yielding the same answer.  We have
    ``conditional consistency''.
  \end{enumerate}
\end{ex}

\begin{ex}
  \label{ex:dempinc}
  Introduce $W$ and $\Omega$ as in \exref{normpivs3}, related by the
  reduced SFM \eqref{wfm}.  We have seen that the marginal fiducial
  distribution for $\Omega$ is the same under the two routes of
  computation.  What about the conditional fiducial distribution of
  $\Sigma$, given $\Omega=\omega$?

  We can reexpress the full  model as
  \begin{eqnarray}
    \label{eq:cc1}
    W &=& \frac{\Omega + E_1}{E_2}\\
    \label{eq:cc2}
    S &=& \Sigma\, E_2.
  \end{eqnarray}

  Again we can identify two routes to construct a conditional
  distribution for $\Sigma$, given $\Omega=\omega$.
  \begin{enumerate}
  \item \label{it:cond3} Form the joint fiducial distribution of
    $(\Omega,\Sigma)$, and condition this on $\Omega=\omega$.

Given data $(x, s)$ (with $\overline x/s = w$) the joint fiducial
    distribution is represented by:
    \begin{eqnarray*}
      \Omega &=& wE_2 - E_1\\
      \Sigma &=& s/E_2.
    \end{eqnarray*}
    So conditioning on $\Omega = \omega$ is equivalent to conditioning on
    \begin{equation}
      \label{eq:bb1}
      wE_2 - E_1 = \omega.
    \end{equation}
    We should therefore condition $E_2$ on this, and then invert
    \eqref{cc2}, so obtaining $\Sigma = s/E_2$, where $E_2$ has its
    distribution conditioned on \eqref{bb1}.
  \item \label{it:cond4} Alternatively we can argue as follows, using
    \eqref{cc1}.  We have observed $W = w$; since we are assuming
    $\Omega = \omega$, we have thus learned
  \begin{equation}
    \label{eq:bb2}
    \frac{\omega + E_1}{E_2} = w.
  \end{equation}
  The conditioning of $E_2$ should therefore be on \eqref{bb2}.
\end{enumerate}
\textcite{demp:incons} showed that we get different answers, depending
on whether we condition on \eqref{bb1} or on \eqref{bb2}.  So here we
have conditioning inconsistency---and it is not clear how we should
resolve it.
\end{ex}
Note that the logical information expressed by \eqref{bb1} and
\eqref{bb2} is the same in both cases.  How then can it matter which
we condition on?  The point is that, when we condition, it is not only
the logical content of the condition that matters, but which partition
of the space it is embedded in.  This is the point of the
``Borel-Kolmogorov paradox'', which concerns conditioning on an event
of probability $0$.  But the paradox can arise even when we have
positive probabilities.

A parable may help.

\begin{ex}
  \label{ex:smith}
  Suppose Mr Smith tells you: ``I have two children, who are not
  twins.''  At this point you regard each of them as equally likely to
  be a boy (B) or a girl (G), independently.  He then says: ``One of
  them is a boy''.  Given this information, what is the probability he
  has two boys?

  \begin{description}
  \item[Argument 1] Initially you assessed 4 equally likely cases: BB,
    BG, GB, GG.  The new information rules out GG, leaving 3 cases,
    just one of which is BB.  The conditional probability is thus 1/3.
  \item[Argument 2] You might consider that, if he had 2 boys, he
    would have said ``They are both boys''.  The fact that he did not
    then implies a conditional probability of 0.
  \end{description}
  Moral: When conditioning on information, we must take account of
  what other information might have been obtained.  Otherwise put, we
  must specify the question (explicit or implicit) that the received
  information answers.  Was it the question ``Do you have a boy?'', or
  the question ``How many boys do you have?''.
\end{ex}

In \exref{dempinc}, the question relevant to \eqref{bb1} is ``What is
the value of $wE_2 - E_1$? (answer: $\omega$).  The question relevant
to \eqref{bb2} is ``What is the value of $(\omega + E_1)/{E_2}$?
(answer: $w$).  Correspondingly the elements of the partition relevant
to \eqref{bb1} are of the form $wE_2 - E_1 = \omega'$, for varying
$\omega'$, while those relevant to \eqref{bb1} are of the form
$(\omega + E_1)/{E_2} = w'$, for varying $w'$.  Only when $w'=w$ and
$\omega' =\omega$ do the answers even contain equivalent logical
information.  Even then, as the partitions differ, so do the
conditional distributions.

In both \exref{normpivs4} and \exref{dempinc}, we wished to condition
on a parameter-function that itself figures in a reduced functional
model.  However, only in \exref{normpivs4} is this model pivotal.
When this is the case, but not more generally. we will obtain the same
partition, and hence the same result, by following each of the two
routes, \itref{cond3} and \itref{cond4}.

\section{Non-simple models}
\label{sec:nonsimp}

\subsection{Ancillary information}
\label{sec:anc}
In many cases there is no simple sufficient statistic.
\textcite{fisher:smsi}, towards the end of the book, suggested---as
usual by means of examples---an alternative approach.  Suppose we can
identify a statistic $S$ that is ancillary, \ie\ has the same
distribution under any $P_\theta$; and a further statistic $T$ such
that, together, $(S,T)$ are equivalent to the full data $X$ (or, more
generally, are jointly sufficient).  Given data $(S,T)=(s,t)$, we can
first restrict attention to the conditional distribution of $T$, given
$S=s$; and then try to identify a pivotal function of $(T,\Theta)$ in
this conditional distribution.  Finally we invert this pivot to obtain
a fiducial distribution.

\begin{ex}
  \label{ex:location}
  Let $X = (X_1,\ldots, X_n)$ arise as a random sample from a general
  location model, with sampling density of the form
  \begin{displaymath}
    f(x \mid \theta) = g(x-\theta).
  \end{displaymath}
  Typically there is no simple sufficient statistic.  However, it
  seems natural to base inference on a location statistic, such as the
  sample mean $\overline X$, or (for $n$ odd) the sample median,
  $\widetilde X$.  But since the sampling distributions of the pivots
  $E_1 = \overline X - \Theta$ and $E_2 = \widetilde X - \Theta$ are
  typically very different, we seem to have a problem of choice.

  This can be resolved as follows.  Let $T = \overline X$, and
  $S = (X_i- \overline X:i=1,\ldots,n)$.  Then $(S,T)$ are together
  equivalent to $X$, $S$ is ancillary, and $E=T-\Theta$ is a pivot,
  both unconditionally and conditionally on $S=s$.  Letting $P^s$
  denote the distribution of $E$ given $S=s$ (the same for all
  $\theta$), a fiducial distribution can be obtained by regarding
  $t-\Theta$ as having distribution $P^s$.  It is easy to show that,
  if we had instead used $T = \widetilde X$ (or any other location
  statistic), the identical fiducial distribution would have been
  obtained.

  The above seemingly well-specified procedure becomes less so when we
  take into account the results of \textcite{basu} that there is,
  typically, a plethora of incommensurate choices of an ancillary to
  condition on.
\end{ex}

\subsection{Non-transitive group models}
\label{sec:nontrans}
We now consider the general case of a group-structured model, as
introduced in \secref{group}, where we do not assume transitivity:
given $x_1,x_2\in\cX$, there may be no $g$ such that $x_2 = gx_1$.
When there is such a $g$ we write $x_1 \approx x_2$.  It is easily
checked that $\approx$ is an equivalence relation on $\cX$: the
equivalence classes under the action of $G$ are termed the {\em
  orbits\/} of $G$ in $\cX$.  Let $S = s(X)$ label orbits.  It is then
readily seen that $S$ is ancillary: this is the {\em group
  ancillary\/}, and, unlike general ancillaries, is essentially
unique.  We choose some arbitrary representative point $x_s$ in the
orbit labelled by $S=s$.  For any $x\in\cX$, there is a unique
$g\in G$ such that $X = g\circ x_s$; we denote this by $t(x)$.  Let
$T=t(X)$.  Then $\Theta^{-1}T$ is a pivot, even conditional on $S=s$.
We can thus apply the construction of \secref{anc} to obtain a
fiducial distribution.  Again, this will coincide with the Bayesian
posterior distribution, based on the right-invariant prior---which may
be easier to compute.

\begin{ex}
  \label{ex:locationbis}
  In \exref{location}, the problem is equivariant under the location
  group, and $S$ and $T$ satisfy the above requirements.  Since the
  right-invariant prior has density element $\d\theta$, the fiducial
  distribution, which could be daunting to compute directly, must have
  density element proportional to
  $\prod_{i=1}^n g(x_i-\theta)\,\d\theta$.
\end{ex}

As noted in \secref{mauldon}, a given problem may be equivariant under
more than one group, and these may induce different fiducial
distributions.  This problem is defined away when we start with a
structural model, which includes specification of the relevant group
$G$.  Then the above recipe yields a unique fiducial distribution.

\subsection{Non-simple functional models}
\label{sec:nonsimpfm}

In a general functional model $X=\Theta\circ E$, we term $x\in\cX$ and
$e\in\cE$ {\em compatible\/} when $x = \theta\circ e$ for some
$\theta\in\cT$.  We first assume {\em invertibility\/}: that such
$\theta$ is unique (compare the group-theoretic concept of exactness),
and write $\theta = x \circ e^{-1}$, noting that in this case $e^{-1}$
is a partial function, operating only on $x$'s that are compatible
with $e$.

Let ${\cal E}_x = \{e:\mbox{$x$ and $e$ are compatible}\}$.  On
observing $X=x$, we learn the logical information $E\in {\cal E}_x$,
but no other logical information about $E$. 

We should thus aim to adjust the distribution of $E$ to account for
this new information, yielding a revised distribution $E \sim P_x$,
say---where $P_x$ is confined to $\cE_x$.  Then
$\Theta = x \circ E^{-1}$ is well-defined, and a fiducial distribution
can be formed by assigning to $E$ the distribution $P_x$ over $\cE_x$.

But how might we compute $P_x$?  As seen in \exref{dempinc},
conditioning on the logical information $E\in \cE_x$ is only
well-defined when $\cE_x$ is embedded in a suitable partition.  We
have to consider what other information we might have obtained, in
other circumstances.  Such information would be of the form
$E\in \cE_y$, as $y$ varies in $\cX$.  Conditioning would thus be
justified when the $\{\cE_y: y\in \cX\}$ form a partition, which will
be the case when, for $x,y\in \cX$, $\cE_x$ and $\cE_y$ are either
identical or disjoint---in which case we term the FM {\em
  partitionable\/}.  Equivalently, there exist, essentially unique,
functions $a(\cdot)$ on $\cX$ , $u(\cdot)$ on $\cE$, such that
$x\in\cX$ and $e\in\cE$ are compatible just when $a(x) = u(e)$.  Since
the observable $X = \Theta\circ E$ is necessarily compatible with the
error variable $E$, $a(X) = u(E)$, and so is ancillary---the {\em
  functional ancillary\/}.  The fiducial distribution of $\Theta$, for
data $X=x$, is now that of $x \circ E^{-1}$, with $E\sim P$
conditioned on $u(E) = a(x)$.

However, a non-partitionable FM does not support unambiguous fiducial
inference.

\subsection{Examples}   
\label{sec:ex3}
We do not have a general necessary and sufficient condition for a
functional model to be partitionable.  This will however hold when the
model is structural, or a reduction of a structural model.

\begin{ex}
  \label{ex:locscale}
  {\bf Location-scale model}\\
  Let ${\cal E} = {\cal X} = \reals^n$,
  $\Theta = (M,\Sigma) \in \cT = \reals \times \reals^+$.  The
  structural model $X = \Theta\circ E$ is given by
  $X_i = M + \Sigma\, E_i$, $i=1,\ldots,n$.  The functional ancillary
  can be taken as $a(x) = ((x_i-\xbar)/s_x: i=1,\ldots,n)$, where
  $s_x^2 = \sum_{i=1}^n (x_i - \xbar)^2/(n-1)$; and $u(e) = a(e)$.
  The fiducial distribution is represented by
  $\Theta = (\xbar - s_x \Ebar/s_e, s_x/s_E)$, where the initial
  distribution of $(\Ebar, s_E)$ is conditioned on
  $(E_i-\Ebar)/s_E = (x_i-\xbar)/s_x$, $i=1,\ldots,n$.  It can
  alternatively be derived as the Bayesian posterior distribution
  based on the right-invariant prior, having density element
  $\d\mu\,\d\sigma/\sigma$.
\end{ex}

\begin{ex}
  \label{ex:locscalered}
  {\bf Reduced structural model}\\
  We have observable $W \in{\cal W} = \{w\in\reals^n:s_w=1\}$,
  parameter $\Omega\in\reals$, error variable
  $E\in{\cal E} = \reals^n$.  The functional model
  $W = \Omega \circ E$ is given by $W_i = (\Omega + E_i)/s_E$,
  $i=1,\ldots,n$.  This is a reduction of the structural model of
  \exref{locscale}, induced by $W = X/s_X$, $\Omega = M/\Sigma$.  It
  is partitionable, with $u(e) = ((e_i-\ebar)/s_e:i=1,\ldots,n)$,
  $a(w) = ((w_i-\wbar):i=1,\ldots,n)$.  The fiducial distribution is
  represented by $\Omega = \wbar s_E - \Ebar$, with the distribution
  of $E$ conditioned on $(E_i-\Ebar)/s_E = w_i-\wbar$, $i=1,\ldots,n$.
  It is not a Bayesian posterior based on the likelihood in the
  reduced model, though it does agree with the marginal for $\Omega$
  in the full fiducial distribution of \exref{locscale} (which is a
  Bayesian posterior).
  
\end{ex}

\begin{ex}
  \label{ex:nonpart}
  {\bf Non-partitionable model}\\
  Let $\cX = \cE = \reals$, $\cT = \reals^+$.  Consider the functional
  model $X=\Theta+E$.  Then ${\cal E}_x = (-\infty, x)$.  Conditioning
  on $E\in{\cal E}_x$ appears, {\em prima facie\/}, straightforward:
  just truncate the initial distribution of $E$ to $(-\infty, x)$.
  However, as the $\{\cE_x: x\in\reals\}$ do not form a partition, it
  is arguable whether this is appropriate.
\end{ex}

\section{Non-invertible models}
\label{sec:nifm}
Consider a functional model $X = \Theta\circ E$, $E\sim P_0$.  Now we
drop the invertibility requirement, so that
$\tau_{x,e} :=\{\theta: x = \theta\circ e\}$ may be a set with more
than one element.

\subsection{Simple non-invertible functional model}
\label{sec:snifm}
We first suppose the model {\em simple\/}, so that any $x\in\cX$ and
$e\in\cE$ are compatible: equivalently, $\tau_{x,e}$ is never empty.

On observing $X=x$, no new logical information is obtained about $E$.
The usual fiducial argument now implies that we can still regard
$E\sim P_0$.  But even were we to know the realised value $e$ of $E$,
we could only infer $\Theta \in \tau_{x,e}$.  In the absence of
knowledge of $e$, the fiducial argument represents our knowledge of
$\Theta$ by $\Theta\in T_x$, where $T_x := \tau_{x,E}$, with
$E\sim P_0$, is a random subset of $\cT$.

This kind of partial probabilistic knowledge, based on random sets,
lies at the heart of the Dempster-Shafer theory of inference
\cite{demp:DScalc}.  Using it, we can go on to define the {\em
  belief\/} and {\em plausibility\/} functions for $\Theta$, after
observing $X=x$:
\begin{eqnarray*}
  \mbox{Bel}_x(\Theta \in A) &=& P_0(T_x \subseteq A)\\
  \mbox{Pl}_x(\Theta \in A) &=& P_0(T_x \cap A \neq \emptyset).
\end{eqnarray*}

\subsection{Recent variations}
\label{sec:recent}
Fiducial theory was largely ignored for many decades.  However recent
years have seen a resurgence of interest, much of it related to
non-invertibility.

\textcite{hannig} carries through an analysis similar to that of
\secref{snifm}, but, in order to finish with a probability
distribution for $\Theta$, adds a further step, in which, given the
compatible set $T_x$, a single value in $T_x$ is selected at random,
from some specified conditional distribution.  There is of course
sensitivity to this specification, and there does not seem to be any
principled way to resolve this.  The theory of inferential models
\cite{martin/liu:book} uses a different auxiliary construction, which
effectively replaces $\mbox{Bel}_x$ by a new belief function
$\mbox{Bel}^*_x$, bounded above by $\mbox{Bel}_x$.  Again there is a
choice of the extra specification.  In both approaches, some guidance
on this may be found by aiming towards compliance with frequentist
(\eg, confidence) properties.

\subsection{General non-invertible functional model}  
\label{sec:gennifm}

We now generalize by allowing $\tau_{x,e} =\emptyset$, equivalent to $x$
and $e$ being incompatible.

On now observing $X=x$, we obtain new logical information about $E$, namely
$$E\in {\cal E}_x := \{e: \tau_{x,e} \neq \emptyset\}.$$
Only when $e\in\cE_x$ could we have made the observation $X=x$ (for
some $\theta\in\cT$).  In order to support fiducial inference, the
initial distribution $P_0$ of $E$ must be adjusted, somehow, to a new
distribution, $P^x$, supported on $\cE_x$.  But how?

Again, things are reasonably straightforward if the model is {\em
  partitionable\/}, \ie, for all $x,x'\in{\cal X}$, ${\cal E}_x$ and
${\cal E}_{x'}$ are either identical or disjoint.  This will this hold
if and only if there exist functions $a(\cdot)$ on $\cX$ and
$u(\cdot)$ on $\cE$, such that $e\in {\cal E}_x$ exactly when
$u(e)=a(x)$.  Then learning $X=x$ is equivalent to learning
$u(E) = a(x)$, and conditioning on this information is unproblematic:
letting $(P_a)$ be the family of conditional distribution of $E$ given
$u(E)=a$ (well-defined under partitionability), we take
$E\sim P_{a(x)}$---which is a distribution supported on $\cE_x$.  We
can finally use this to define the distribution of the random set
$T_{x}$ (and so $\mbox{Bel}_x$, $\mbox{Pl}_x$).  However, when the
model is not partitionable it is not clear how (or indeed whether) to
construct $P^x$.  \textcite{hannig} suggests ways of identifying a
suitable function on which to condition, but there typically remains a
multiplicity of apparently reasonable choices.

\begin{ex}
  \label{ex:hannig}
  Consider the following functional model (a variation on Example~4 of
  \textcite[Examle~4]{hannig}):
  \begin{eqnarray}
    \label{eq:han1}
    X_1 &=& \Theta_1/E_1\\
    \label{eq:han2}
    X_2 &=& (\Theta_1 + \Theta_2)/E_2\\
    \label{eq:han3}
    X_3 &=& (\Theta_1 + 2\Theta_2)/E_3.
  \end{eqnarray}
  Fixing data $X=x$ and solving \eqref{han1} and \eqref{han2}, we
  obtain
  \begin{eqnarray}
    \label{eq:han4}
    \Theta_1 &=& x_1E_1\\
    \label{eq:han5}
    \Theta_2 &=& x_2 E_2 -x_1 E_1,
  \end{eqnarray}
  and then inserting these in \eqref{han3} we get the compatibility
  condition, $E \in \cE_x$, expressed as
  \begin{equation}
    \label{eq:han6}
    \frac{2x_2E_2 - x_1 E_1}{E_3}=  x_3.
  \end{equation}
  A similar analysis that starts by solving \eqref{han2} and
  \eqref{han3} yields
  \begin{equation}
    \label{eq:han7}
    \frac{2x_2E_2 - x_3 E_3}{E_1}=  x_1,
  \end{equation}
  Both \eqref{han6} and \eqref{han7} are (necessarily) equivalent to
  each other, and to
  \begin{equation}
    \label{eq:han8}
    x_1E_1 - 2x_2  E_2 + x_3 E_3 = 0.
  \end{equation}
  The partitions generated by the variables on the left-hand sides of
  \eqref{han6}, \eqref{han7} and \eqref{han8} are all different, so
  that conditioning \eqref{han4} and \eqref{han5} on them will give
  different answers.  And indeed, since the model is
  non-partitionable, there is no correct answer as to how (or whether)
  we should condition.  To see non-partitionability directly, note
  that, as expressed by \eqref{han8}, each $\cE_x$ is a plane in
  $\reals^3$.  When not identical, any two such planes must intersect
  in a line, so can not be disjoint.
  
  \textcite{hannig} notes that, in examples such as this, it matters
  how we condition on the information $E \in \cE_x$, and makes some
  {\em ad hoc\/} recommendations.  But in view of
  non-partitionability, it could be argued that no conditioning of any
  kind is justifiable, and that the model simply does not support
  fiducial inference.
\end{ex}

\begin{ex}
\label{ex:nifm}
Let $\cX = \cE = [0,1]^n$, $\cT = [0,1]$.  Under $P_0$,
$(E_i:i=1,\ldots,n)$ are independently uniform over $[0,1]$.  The
functional model is given by $X_i = \bbl(E_i \leq \Theta)$,
$i=1,\ldots,n$.  Then, when $\Theta=\theta$, the $X_i$ are $n$
independent Bernoulli($\theta$) variables.  This functional model is
the basis of Example~6 of \textcite{hannig}.

We see that $e$ and $x$ are compatible ($\tau_{x,e}\neq\emptyset$)
when $x_i=1, x_j = 0$ if and only if $e_i < e_j$.  Thus on observing
$X=x$, we learn
\begin{equation}
  \label{eq:ee}  
  E \in \cE_{x}:= \{e: e_i < e_j \mbox{ just when }x_i=1, x_j = 0\}, 
\end{equation}
and then $\Theta$ lies in the random interval between two order statistics:
$$T_{E} := [\,E_{(r)}, E_{(r+1)}\,)\quad (r=\sum x_i).$$

We have to confine the distribution of $E$ to $\cE_x$, but how?
Noting that $P_0(\cE_x)>0$, an obvious approach is simply to truncate
$P_0$ to the set $\cE_x$.  We may then note that \eqref{ee} is a
condition on the way in which the $(E_i)$ are ordered, which is
independent of their order-statistic, which is what determines $T_E$.
So the distribution of the random interval $T_E$ will be unaffected by
the truncation to \eqref{ee}---allowing us to use its unconditional
distribution, based on the order statistics of a random sample from
the uniform distribution on $[0,1]$.  This is the approach of
\textcite{hannig}.

Nevertheless, this model is not partitionable, as may be seen by
noting that, when $x = {\bf 1}$ (the vector with all $x_i=1$), we have
${\cE}_{\bf 1} = \cE$---which is not disjoint from or identical with
any other $\cE_{x}$.  So, in the light of examples such as
\exref{smith}, and in the absence of a clear question that is answered
by the information $E\in\cE_x$, it is debatable whether the above
argument is appropriate.  Once again there seems to be no fully
justifiable fiducial inference available.

\end{ex}

\section{Concluding comments}  
\label{sec:conc}
While the fiducial argument has some {\em prima facie\/} appeal, all
attempts to formulate a fully coherent theory of fiducial inference
have fallen foul of inconsistency and counter-examples.  The
investigations of \textcite{apd/ms:funcmod}, based on functional
models, were an attempt to see just how far the theory could be taken
before crashing onto the rocks---but crash it eventually did.  Many of
the difficulties are associated with the need to specify,
unambiguously, a relevant partition for performing probabilistic
conditioning.  In some cases, as in \exref{dempinc}, two equally
natural routes to take account of new or assumed information lead to
different partitions and hence conflicting fiducial distributions---a
parallel inconsistency in predictive inference was exhibited by
\textcite{apd/jw}.  In other cases, as in \secref{gennifm}, there is
no natural embedding of the information obtained within any partition
whatsoever, rendering fiducial inference undefined.

Nevertheless, even though methods derived from fiducial theory may
have limited validity from a fully principled theoretical standpoint,
that is not to deny that they may prove useful for other
purposes---for example, for constructing exact or approximate
confidence regions.  But---simple pivotal cases apart---there is no
guarantee that this will be the case, so that further investigations
are required in individual cases.


\begin{thebibliography}{}

\bibitem[\protect\citeAY{Angus}{1994}]{angus}
Angus, J.~E. (1994).
\newblock The probability integral transform and related results.
\newblock {\em SIAM Review}, {\bf 36}, 652--4.

\bibitem[\protect\citeAY{Barnard}{1980}]{barnard}
Barnard, G.~A. (1980).
\newblock Pivotal inference and the {B}ayesian controversy.
\newblock {\em Trabajos de Estadistica y de Investigacion Operativa}, {\bf 31},
  295--318.

\bibitem[\protect\citeAY{Basu}{1959}]{basu}
Basu, D. (1959).
\newblock The family of ancillary statistics.
\newblock {\em Sankhy\=a}, {\bf 21}, 247--56.

\bibitem[\protect\citeAY{Bondar}{1972}]{bondar}
Bondar, J.~V. (1972).
\newblock Structural distributions without exact transitivity.
\newblock {\em Annals of Mathematical Statistics}, {\bf 43}, 326--39.

\bibitem[\protect\citeAY{Buehler}{1983}]{rjb:ess}
Buehler, R.~J. (1983).
\newblock Fiducial inference.
\newblock In {\em Encyclopedia of Statistical Sciences, Volume~3}, (ed.\
  S.~Kotz, N.~L. Johnson, and C.~B. Read), pp.~76--81. Wiley-Interscience.
\newblock \\\href{https://doi.org/10.1002/9781118445112.stat01529}{\tt
  DOI:10.1002/9781118445112.stat01529}.

\bibitem[\protect\citeAY{Dawid and Stone}{1982}]{apd/ms:funcmod}
Dawid, A.~P. and Stone, M. (1982).
\newblock The functional-model basis of fiducial inference (with {Discussion)}.
\newblock {\em Annals of Statistics}, {\bf 10}, 1054--74.

\bibitem[\protect\citeAY{Dawid {\it et~al}.}{1973}]{dsz}
Dawid, A.~P., Stone, M., and Zidek, J.~V. (1973).
\newblock Marginalization paradoxes in {B}ayesian and structural inference
  (with {D}iscussion).
\newblock {\em Journal of the Royal Statistical Society. Series B
  (Methodological)}, {\bf 35}, 189--233.

\bibitem[\protect\citeAY{Dawid and Wang}{1993}]{apd/jw}
Dawid, A.~P. and Wang, J. (1993).
\newblock Fiducial prediction and semi-{B}ayesian inference.
\newblock {\em Annals of Statistics}, {\bf 21}, 1119--38.

\bibitem[\protect\citeAY{Dempster}{1963}]{demp:incons}
Dempster, A.~P. (1963).
\newblock Further examples of inconsistencies in the fiducial argument.
\newblock {\em Annals of Mathematical Statistics}, {\bf 34}, 884--91.

\bibitem[\protect\citeAY{Dempster}{2008}]{demp:DScalc}
Dempster, A.~P. (2008).
\newblock The {Dempster-Shafer} calculus for statisticians.
\newblock {\em International Journal of Approximate Reasoning}, {\bf 48},
  365--77.

\bibitem[\protect\citeAY{Edwards}{1983}]{awfe:ess}
Edwards, A. W.~F. (1983).
\newblock Fiducial distributions.
\newblock In {\em Encyclopedia of Statistical Sciences, Volume~3}, (ed.\
  S.~Kotz, N.~L. Johnson, and C.~B. Read), pp.~70--6. Wiley-Interscience.
\newblock \\\href{https://doi.org/10.1002/0471667196.ess0775.pub2}{\tt
  DOI:10.1002/0471667196.ess0775.pub2}.

\bibitem[\protect\citeAY{Fisher}{1915}]{fisher:cor}
Fisher, R.~A. (1915).
\newblock Frequency distribution of the values of the correlation coefficient
  in samples from an indefinitely large population.
\newblock {\em Biometrika}, {\bf 10}, 507--21.

\bibitem[\protect\citeAY{Fisher}{1930}]{fisher:invprob}
Fisher, R.~A. (1930).
\newblock Inverse probability.
\newblock {\em Mathematical Proceediings of the ~Cambridge Philosophical
  Society}, {\bf 26}, 528--35.

\bibitem[\protect\citeAY{Fisher}{1956}]{fisher:smsi}
Fisher, R.~A. (1956).
\newblock {\em Statistical Methods and Scientific Inference}. Oliver and Boyd,
  Edinburgh.
\newblock Third Edition, Hafner, New York, 1973.

\bibitem[\protect\citeAY{Fraser}{1961}]{fraser:1961}
Fraser, D. A.~S. (1961).
\newblock On fiducial inference.
\newblock {\em Annals of Mathematical Statistics}, {\bf 32}, 661--76.

\bibitem[\protect\citeAY{Fraser}{1968}]{fraser:book}
Fraser, D. A.~S. (1968).
\newblock {\em The Structure of Inference}. Wiley, New York.

\bibitem[\protect\citeAY{Hacking}{1965}]{hacking:book}
Hacking, I. (1965).
\newblock {\em Logic of Statistical Inference}. Cambridge University Press,
  Cambridge.

\bibitem[\protect\citeAY{Hannig}{2009}]{hannig}
Hannig, J. (2009).
\newblock On generalized fiducial inference.
\newblock {\em Statistica Sinica}, {\bf 19}, 491--544.

\bibitem[\protect\citeAY{Harris and Harding}{1984}]{harris}
Harris, R.~R. and Harding, E.~F. (1984).
\newblock The fiducial argument and {H}acking's principle of irrelevance.
\newblock {\em Journal of Applied Statistics}, {\bf 11}, 170--81.

\bibitem[\protect\citeAY{Lindley}{1958}]{dvl:fidbayes}
Lindley, D.~V. (1958).
\newblock Fiducial distributions and {B}ayes' theorem.
\newblock {\em Journal of the Royal Statistical Society, Series~B}, {\bf 20},
  102--7.

\bibitem[\protect\citeAY{Martin and Liu}{2016}]{martin/liu:book}
Martin, R. and Liu, C. (2016).
\newblock {\em Inferential Models}. Chapman and Hall/CRC, New York.

\bibitem[\protect\citeAY{Mauldon}{1955}]{mauldon}
Mauldon, J.~G. (1955).
\newblock Pivotal quantities for {W}ishart's and related distributions, and a
  paradox in fiducial theory.
\newblock {\em Journal of the Royal Statistical Society, Series~B}, {\bf 17},
  79--85.

\bibitem[\protect\citeAY{Savage}{1961}]{savage:eggs}
Savage, L.~J. (1961).
\newblock The foundations of statistics reconsidered.
\newblock In {\em Proceedings of the Fourth Berkeley Symposium on Mathematical
  Statistics and Probability, Volume 1}, pp.~575--86. University of California
  Press.

\bibitem[\protect\citeAY{Stone}{1983}]{ms:ess}
Stone, M. (1983).
\newblock Fiducial probability.
\newblock In {\em Encyclopedia of Statistical Sciences, Volume~3}, (ed.\
  S.~Kotz, N.~L. Johnson, and C.~B. Read), pp.~81--6. Wiley-Interscience.
\newblock \\\href{https://doi.org/10.1002/9781118445112.stat01530}{\tt
  DOI:10.1002/9781118445112.stat01530}.

\bibitem[\protect\citeAY{Wallace}{1980}]{wallace:bffc}
Wallace, D.~L. (1980).
\newblock The {Behrens-Fisher} and {Fieller-Creasy} problems.
\newblock In {\em R. A. Fisher: An Appreciation}, (ed.\ S.~E. Fienberg and
  D.~V. Hinkley).
\newblock Springer New York, New York, NY.

\bibitem[\protect\citeAY{Wilkinson}{1977}]{wilkinson}
Wilkinson, G.~N. (1977).
\newblock On resolving the controversy in statistical inference (with
  {D}iscussion).
\newblock {\em Journal of the Royal Statistical Society, Series~B}, {\bf 39},
  119--71.

\bibitem[\protect\citeAY{Zabell}{1992}]{zabell}
Zabell, S.~L. (1992).
\newblock {R. A. F}isher and the fiducial argument.
\newblock {\em Statistical Science}, {\bf 7}, 369--87.

\end{thebibliography}
\end{document}